\catcode`\@=11

\font\tenmsa=msam10 \font\sevenmsa=msam7 \font\fivemsa=msam5
\font\tenmsb=msbm10
\font\sevenmsb=msbm7 \font\fivemsb=msbm5 \newfam\msafam \newfam\msbfam
\textfont\msafam=\tenmsa \scriptfont\msafam=\sevenmsa
\scriptscriptfont\msafam=\fivemsa \textfont\msbfam=\tenmsb
\scriptfont\msbfam=\sevenmsb \scriptscriptfont\msbfam=\fivemsb

\def\hexnumber@#1{\ifnum#1<10 \number#1\else \ifnum#1=10 A\else\ifnum#1=11
 B\else\ifnum#1=12 C\else \ifnum#1=13 D\else\ifnum#1=14 E\else\ifnum#1=15
 F\fi\fi\fi\fi\fi\fi\fi}

\def\msa@{\hexnumber@\msafam} \def\msb@{\hexnumber@\msbfam}
\mathchardef\boxdot="2\msa@00 \mathchardef\boxplus="2\msa@01
\mathchardef\boxtimes="2\msa@02 \mathchardef\square="0\msa@03
\mathchardef\blacksquare="0\msa@04 \mathchardef\centerdot="2\msa@05
\mathchardef\lozenge="0\msa@06 \mathchardef\blacklozenge="0\msa@07
\mathchardef\circlearrowright="3\msa@08 \mathchardef\circlearrowleft="3\msa@09
\mathchardef\rightleftharpoons="3\msa@0A
\mathchardef\leftrightharpoons="3\msa@0B \mathchardef\boxminus="2\msa@0C
\mathchardef\Vdash="3\msa@0D \mathchardef\Vvdash="3\msa@0E
\mathchardef\vDash="3\msa@0F \mathchardef\twoheadrightarrow="3\msa@10
\mathchardef\twoheadleftarrow="3\msa@11 \mathchardef\leftleftarrows="3\msa@12
\mathchardef\rightrightarrows="3\msa@13 \mathchardef\upuparrows="3\msa@14
\mathchardef\downdownarrows="3\msa@15 \mathchardef\upharpoonright="3\msa@16
 \mathchardef\downharpoonright="3\msa@17
\mathchardef\upharpoonleft="3\msa@18 \mathchardef\downharpoonleft="3\msa@19
\mathchardef\rightarrowtail="3\msa@1A \mathchardef\leftarrowtail="3\msa@1B
\mathchardef\leftrightarrows="3\msa@1C \mathchardef\rightleftarrows="3\msa@1D
\mathchardef\Lsh="3\msa@1E \mathchardef\Rsh="3\msa@1F
\mathchardef\rightsquigarrow="3\msa@20
\mathchardef\leftrightsquigarrow="3\msa@21 \mathchardef\looparrowleft="3\msa@22
\mathchardef\looparrowright="3\msa@23 \mathchardef\circeq="3\msa@24
\mathchardef\succsim="3\msa@25 \mathchardef\gtrsim="3\msa@26
\mathchardef\gtrapprox="3\msa@27 \mathchardef\multimap="3\msa@28
\mathchardef\therefore="3\msa@29 \mathchardef\because="3\msa@2A
\mathchardef\doteqdot="3\msa@2B 
\mathchardef\traceiangleq="3\msa@2C \mathchardef\precsim="3\msa@2D
\mathchardef\lesssim="3\msa@2E \mathchardef\lessapprox="3\msa@2F
\mathchardef\eqslantless="3\msa@30 \mathchardef\eqslantgtr="3\msa@31
\mathchardef\curlyeqprec="3\msa@32 \mathchardef\curlyeqsucc="3\msa@33
\mathchardef\preccurlyeq="3\msa@34 \mathchardef\leqq="3\msa@35
\mathchardef\leqslant="3\msa@36 \mathchardef\lessgtr="3\msa@37
\mathchardef\backprime="0\msa@38 \mathchardef\risingdotseq="3\msa@3A
\mathchardef\fallingdotseq="3\msa@3B \mathchardef\succcurlyeq="3\msa@3C
\mathchardef\geqq="3\msa@3D \mathchardef\geqslant="3\msa@3E
\mathchardef\gtrless="3\msa@3F \mathchardef\sqsubset="3\msa@40
\mathchardef\sqsupset="3\msa@41
\mathchardef\trianglelefteq="3\msa@45 \mathchardef\bigstar="0\msa@46
\mathchardef\between="3\msa@47 \mathchardef\blacktriangledown="0\msa@48
\mathchardef\blacktriangleright="3\msa@49
\mathchardef\blacktriangleleft="3\msa@4A
\mathchardef\blacktriangle="0\msa@4E \mathchardef\triangledown="0\msa@4F
\mathchardef\eqcirc="3\msa@50 \mathchardef\lesseqgtr="3\msa@51
\mathchardef\gtreqless="3\msa@52 \mathchardef\lesseqqgtr="3\msa@53
\mathchardef\gtreqqless="3\msa@54 \mathchardef\Rrightarrow="3\msa@56
\mathchardef\Lleftarrow="3\msa@57 \mathchardef\veebar="2\msa@59
\mathchardef\barwedge="2\msa@5A \mathchardef\doublebarwedge="2\msa@5B
\mathchardef\angle="0\msa@5C \mathchardef\measuredangle="0\msa@5D
\mathchardef\sphericalangle="0\msa@5E \mathchardef\varpropto="3\msa@5F
\mathchardef\smallsmile="3\msa@60 \mathchardef\smallfrown="3\msa@61
\mathchardef\Subset="3\msa@62 \mathchardef\Supset="3\msa@63
\mathchardef\Cup="2\msa@64  \mathchardef\Cap="2\msa@65
 \mathchardef\curlywedge="2\msa@66
\mathchardef\curlyvee="2\msa@67 \mathchardef\leftthreetimes="2\msa@68
\mathchardef\rightthreetimes="2\msa@69 \mathchardef\subseteqq="3\msa@6A
\mathchardef\supseteqq="3\msa@6B \mathchardef\bumpeq="3\msa@6C
\mathchardef\Bumpeq="3\msa@6D \mathchardef\lll="3\msa@6E 
\mathchardef\ggg="3\msa@6F  \mathchardef\circledS="0\msa@73
\mathchardef\pitchfork="3\msa@74 \mathchardef\dotplus="2\msa@75
\mathchardef\backsim="3\msa@76 \mathchardef\backsimeq="3\msa@77
\mathchardef\complement="0\msa@7B \mathchardef\intercal="2\msa@7C
\mathchardef\circledcirc="2\msa@7D \mathchardef\circledast="2\msa@7E
\mathchardef\circleddash="2\msa@7F \def\ulcorner{\delimiter"4\msa@70\msa@70 }
\def\urcorner{\delimiter"5\msa@71\msa@71 }
\def\llcorner{\delimiter"4\msa@78\msa@78 }
\def\lrcorner{\delimiter"5\msa@79\msa@79 } \def\yen{\mathhexbox\msa@55 }
\def\checkmark{\mathhexbox\msa@58 } \def\circledR{\mathhexbox\msa@72 }
\def\maltese{\mathhexbox\msa@7A } \mathchardef\lvertneqq="3\msb@00
\mathchardef\gvertneqq="3\msb@01 \mathchardef\nleq="3\msb@02
\mathchardef\ngeq="3\msb@03 \mathchardef\nless="3\msb@04
\mathchardef\ngtr="3\msb@05 \mathchardef\nprec="3\msb@06
\mathchardef\nsucc="3\msb@07 \mathchardef\lneqq="3\msb@08
\mathchardef\gneqq="3\msb@09 \mathchardef\nleqslant="3\msb@0A
\mathchardef\ngeqslant="3\msb@0B \mathchardef\lneq="3\msb@0C
\mathchardef\gneq="3\msb@0D \mathchardef\npreceq="3\msb@0E
\mathchardef\nsucceq="3\msb@0F \mathchardef\precnsim="3\msb@10
\mathchardef\succnsim="3\msb@11 \mathchardef\lnsim="3\msb@12
\mathchardef\gnsim="3\msb@13 \mathchardef\nleqq="3\msb@14
\mathchardef\ngeqq="3\msb@15 \mathchardef\precneqq="3\msb@16
\mathchardef\succneqq="3\msb@17 \mathchardef\precnapprox="3\msb@18
\mathchardef\succnapprox="3\msb@19 \mathchardef\lnapprox="3\msb@1A
\mathchardef\gnapprox="3\msb@1B \mathchardef\nsim="3\msb@1C
\mathchardef\napprox="3\msb@1D
\mathchardef\nsupseteqq="3\msb@23 \mathchardef\subsetneqq="3\msb@24
\mathchardef\supsetneqq="3\msb@25
\mathchardef\supsetneq="3\msb@29 \mathchardef\nsubseteq="3\msb@2A
\mathchardef\nsupseteq="3\msb@2B \mathchardef\nparallel="3\msb@2C
\mathchardef\nmid="3\msb@2D \mathchardef\nshortmid="3\msb@2E
\mathchardef\nshortparallel="3\msb@2F \mathchardef\nvdash="3\msb@30
\mathchardef\nVdash="3\msb@31 \mathchardef\nvDash="3\msb@32
\mathchardef\nVDash="3\msb@33 \mathchardef\ntrianglerighteq="3\msb@34
\mathchardef\ntrianglelefteq="3\msb@35 \mathchardef\ntriangleleft="3\msb@36
\mathchardef\ntriangleright="3\msb@37 \mathchardef\nleftarrow="3\msb@38
\mathchardef\nrightarrow="3\msb@39 \mathchardef\nLeftarrow="3\msb@3A
\mathchardef\nRightarrow="3\msb@3B \mathchardef\nLeftrightarrow="3\msb@3C
\mathchardef\nleftrightarrow="3\msb@3D \mathchardef\divideontimes="2\msb@3E
\mathchardef\varnothing="0\msb@3F \mathchardef\nexists="0\msb@40
\mathchardef\mho="0\msb@66 \mathchardef\thorn="0\msb@67
\mathchardef\beth="0\msb@69 \mathchardef\gimel="0\msb@6A
\mathchardef\daleth="0\msb@6B \mathchardef\lessdot="3\msb@6C
\mathchardef\gtrdot="3\msb@6D \mathchardef\ltimes="2\msb@6E
\mathchardef\rtimes="2\msb@6F \mathchardef\shortmid="3\msb@70
\mathchardef\shortparallel="3\msb@71 \mathchardef\smallsetminus="2\msb@72
\mathchardef\thicksim="3\msb@73 \mathchardef\thickapprox="3\msb@74
\mathchardef\approxeq="3\msb@75 \mathchardef\succapprox="3\msb@76
\mathchardef\precapprox="3\msb@77 \mathchardef\curvearrowleft="3\msb@78
\mathchardef\curvearrowright="3\msb@79 \mathchardef\digamma="0\msb@7A
\mathchardef\varkappa="0\msb@7B \mathchardef\hslash="0\msb@7D
\mathchardef\hbar="0\msb@7E \mathchardef\backepsilon="3\msb@7F
\def\Bbb{\ifmmode\let\next\Bbb@\else
\def\next{\errmessage{Use \string\Bbb\space only in math mode}}\fi\next}
\def\Bbb@#1{{\Bbb@@{#1}}} \def\Bbb@@#1{\fam\msbfam#1}

\catcode`\@=\active

 \def\CC{\hbox{{$\cal C$}}}

\def\CR{\hbox{{$\cal R$}}} 
 \def\CV{\hbox{{$\cal V$}}}
\def\CM{\hbox{{$\cal M$}}}

\def\lform{\hbox{$\sqcup$}\llap{\hbox{$\sqcap$}}}

\def\Z{{\Bbb Z}}

\def\eps{{\epsilon}}

\def\aut{{\rm Aut\, }}

\def\lcross{{>\!\!\!\triangleleft}}
\def\rcocross{{\blacktriangleright\!\!<}}
\def\lcocross{{>\!\!\blacktriangleleft}}

\def\rbiprod{{\cdot\kern-.33em\triangleright\!\!\!<}}
\def\lbiprod{{>\!\!\!\triangleleft\kern-.33em\cdot\, }}

\def\tens{\mathop{\otimes}}

\def\la{{\triangleright}}\def\ra{{\triangleleft}}

\def\isom{{\cong}}

\def\Nat{{\rm Nat}}

\def\<{\langle}
\def\>{\rangle}

\def\equad{\kern -1.7em}

\def\eqn#1#2{\begin{equation}#2\label{#1}\end{equation}}

\def\o{{}_{\scriptscriptstyle(1)}}
\def\t{{}_{\scriptscriptstyle(2)}}
\def\th{{}_{\scriptscriptstyle(3)}}
\def\fo{{}_{\scriptscriptstyle(4)}}

\def\bo{{}^{\bar{\scriptscriptstyle(1)}}}
\def\bt{{}^{\bar{\scriptscriptstyle(2)}}}

\def\und#1{{\underline {#1}}}

\def\uo{{{}^{\scriptscriptstyle(1)}}}
\def\ut{{{}^{\scriptscriptstyle(2)}}}

\def\Bo{{{}_{\und{\scriptscriptstyle(1)}}}}
\def\Bt{{{}_{\und{\scriptscriptstyle(2)}}}}

\def\new#1{\goodbreak\goodbreak\bigskip
\noindent{\bf #1}}

\def\text#1{\mbox{\rm #1}}
\def\note#1{}

\def\blacksquare{{\lform}}
\def\frac#1#2{{{#1\over#2}}}

\def\proof{\goodbreak\noindent{\bf Proof\quad}}

\def\endproof{{\ $\lform$}\bigskip }

\def\alignn#1#2{\begin{eqnarray}\label{#1}#2
\end{eqnarray}}
\def\cmath#1{\[\begin{array}{c} #1 \end{array}\]}
\def\ceqn#1#2{\begin{equation}\label{#1}
\begin{array}{c}#2\end{array}\end{equation}}

\def\vecu{{\bf u}}

\documentstyle[11pt,epsf]{article} \textheight 23.6cm \textwidth 16cm
\evensidemargin
0in \topskip 28pt

\newtheorem{lemma}{Lemma} \newtheorem{propos}[lemma]{Proposition}
\newtheorem{example}[lemma]{Example}

\begin{document}\baselineskip 20pt

{\ }\hskip 4.7in  DAMTP/95-19
\vspace{.2in}

\begin{center} {\LARGE  COMMENTS ON BOSONISATION AND BIPRODUCTS} \\
\baselineskip 13pt
{\ } {\ }\\ S.  Majid\footnote{Royal Society University Research Fellow and
Fellow of Pembroke College, Cambridge}\\{\ }\\ Department of Applied
Mathematics \& Theoretical Physics\\ University of Cambridge, Cambridge CB3
9EW, UK\\
+\\
Department of Mathematics, Harvard University\\
1 Oxford Street, Cambridge, MA02138, USA\footnote{During the calendar years
1995+1996}
\end{center}
\begin{center}
April 9, 1994
\end{center}
\begin{quote} \baselineskip 13pt
\noindent{\bf Abstract}  We clarify the relation  between the `bosonisation'
construction (due to the author)  which can be used to turn a Hopf algebra $B$
in ${}_H\CM$ or ${}^H\CM$ into an equivalent ordinary Hopf algebra, and a
version of Radford's theorem (also due in this form to the author) which does
the same for $B$ in ${}^H_H\CM$. We also comment on reconstruction from the
category of $B$-comodules.
\end{quote}

Recently there appeared an interesting paper by Fischman and
Montgomery\cite{FisMon:sch} in which some abstract constructions due
to the author for Hopf algebras in braided categories
\cite{Ma:bg}\cite{Ma:dou}\cite{Ma:skl}\cite{Ma:bos} were applied
to obtain a Schur's double centraliser theorem for color Lie algebras.
Unfortunately, in the introduction and preliminaries of their paper, some
incorrect statements are made which could cause confusion and which we would
like to correct here in this note. The most serious is the identification
${}^H\CM={}^H_H\CM$ which we cover in Section~8. The main results about
color-Lie algebras, etc., in \cite{FisMon:sch} are, fortunately, not affected.
Finally, we add some related remarks about a recent preprint of Pareigis
\cite{Par:rec}.

Most of these comments were made to Susan Montgomery at the LMS meeting on
Noncommutative Rings, Durham, UK in July 1992, where we were shown the first
draft of \cite{FisMon:sch} and pointed out
\cite{Ma:bg}\cite{Ma:dou}\cite{Ma:skl}\cite{Ma:bos}. We work over a ground
field $k$ and use the usual notations $S,\eps$ for the  antipode and counit,
and the notation $\Delta h=\sum h\o\tens h\t$ for the coproduct applied to an
element $h\in H$\cite{Swe:hop}. We use the symbols $\lcross$, etc., for smash
products, $\lcocross$ for smash coproducts and $\lbiprod$ when both are made
simultaneously.

\new{1. \sl Dual quasitriangular structures.} Let $H$ be a dual quasitriangular
(or coquasitriangular) Hopf algebra, i.e. equipped with a
convolution-invertible linear map $\CR:H\tens H\to k$ obeying
\ceqn{dqua}{ \CR(h\tens gf)=\sum \CR(h\o\tens f)\CR(h\t,g),\quad \CR(fg\tens
h)=\sum \CR(f\tens h\o)\CR(g\tens h\t)\\
\sum g\o h\o\CR(h\t\tens g\t)=\sum \CR(h\o\tens g\o) h\t g\t}
These axioms are nothing other than the axioms of a quasitriangular structure
introduced by Drinfeld\cite{Dri}, with product replaced by convolution product
of $\CR$ as a linear functional. They were formulated by the author in
1989\cite[below~Thm~2.1]{Ma:rep},\cite[Thm~4.1]{Ma:pro}
cf\cite[Sec.~7.4-7.4]{Ma:qua} as equivalent under Tannaka-Krein reconstruction
to the category $\CM^H$ of right comodules becoming  braided in the sense of
\cite{JoyStr:bra}. The result was generalised in \cite{Ma:tan} to the dual
quasitriangular quasi-Hopf algebra setting in the Fall of 1990. One can use
left comodules ${}^H\CM$ just as well.

The basic properties of dual quasitriangular Hopf algebras as algebraic objects
are due to the author in \cite[Appendix]{Ma:bg}, again from the Fall of 1990.
For example

\begin{propos}\cite[Prop~A.5]{Ma:bg} Let $H$ be a dual quasitriangular Hopf
algebra. Then the square of the antipode is inner in the convolution algebra
$H\to H$ and hence the antipode is bijective.
\end{propos}

This is nothing more than the dual of standard results (due ultimately to V.G.
Drinfeld\cite{Dri:alm}) for quasitriangular Hopf algebras. From such basic
results in \cite{Ma:bg} we see that the assumption that the antipode of $H$ is
bijective in the paper of Fischman and Montgomery (e.g. in
\cite[Thm~2.15]{FisMon:sch}) should be deleted as superfluous.

The matrix bialgebras $A(R)$ of \cite{FRT:lie} cf.\cite{Dri} associated to a
matrix $R$ are dual quasitriangular whenever $R$ obeys the so-called quantum
Yang-Baxter equations or Artin braid relations. This result is due to the
author in 1989 in \cite[p. 141]{Ma:seq} as  $\CR:A(R)\to A(R)^*$ and was
elaborated further as $\CR:A(R)\tens A(R)\to k$ by Larson and Towber
\cite{LarTow:two} under another name.

\new{2. \sl Braided groups.} In case the attribution in \cite[p.
596]{FisMon:sch} is not clear, we would like to stress that the theory of
`braided groups' or bialgebras and Hopf algebras $B$ living in a braided
category is due to the present author in 1989 and V. Lyubashenko in 1990, and
arose in both cases in connection with conformal field theory and generalised
Tannaka-Krein reconstruction theorems. See \cite{LyuMa:bra} and
\cite{Ma:introm} for a review. The theory is more subtle than the theory of
Hopf algebras in symmetric categories, which have an older history such as
\cite{Par:non}.
If $A,B$ are two algebras in the braided category, one checks first that
$A\tens B$ has a {\em braided tensor product algebra} structure
\cite{Ma:eul}\cite{Ma:bg} with product morphism
$(\cdot_A\tens\cdot_B)\circ\Psi_{B,A}$ where $\Psi_{B,A}:B\tens A\to A\tens B$
is the braiding morphism between any two objects in the category. Note that
objects in a category need not be sets with elements. A bialgebra in a braided
category is
an algebra $B$ which is also a coalgebra with coproduct and counit
$\Delta_B,\eps_B$ algebra homomorphisms. Here $\Delta_B: B\to B\tens B$ to the
braided tensor product algebra above. The morphisms $B\to B$ form an algebra
under convolution and an antipode is defined as usual as an inverse for the
identity. Basic lemmas about Hopf algebras in braided categories or {\em
braided groups} are due to the author. For example,

\begin{propos}\cite[Fig.~2]{Ma:tra} The antipode $S_B$ of a braided group $B$
is braided-antimultiplicative in the sense
\[ S_B\circ \cdot_B=\cdot_B\circ\Psi_{B,B}\circ(S_B\tens S_B),\quad
\Delta_B\circ
S_B=(S_B\tens S_B)\circ\Psi_{B,B}\circ\Delta_B.\]
\end{propos}

The proof\cite{Ma:tra} in the braided case is non trivial and makes use of the
coherence theorem for braided categories\cite{JoyStr:bra} in the form of a
diagrammatic notation in which algebraic information `flows' along braids and
tangles. Further basic properties such tensor products of modules and comodules
in the category, smash products and coproducts in the category, etc, can be
found in \cite{Ma:tra}\cite{Ma:bos} where we  give diagrammatic proofs for the
module versions. The comodule versions are the same with all diagrams turned
up-side-down. There is also a theory of quasitriangular structures and dual
quasitriangular structures for such categorical braided groups.

\new{3. \sl Braided groups in $\CM^H$.} In case it is not clear from
\cite[Remark~1.6(e)]{FisMon:sch}, the theory and construction of bialgebras and
Hopf algebras in the specific braided category of comodules of a dual
quasitriangular Hopf $H$ is due to the author in \cite{Ma:eul}\cite{Ma:bg}. It
was presented at St Petersburg in September, 1990 and at the AMS Special
Session on Hopf Algebras, San Francisco, January 1991, where it was presented
to the `Hopf algebras community' including Susan Montgomery. This was new  even
for the symmetric (unbraided) case of primary interest in \cite{FisMon:sch}. We
use the notation $\Delta_B b=\sum b\Bo\tens b\Bt$ for the coproduct of $B$ and
$\beta(b)=\sum b\bo\tens b\bt$ for a left or right coaction on it.

Basic properties are now $\Delta_B$ an algebra homomorphism to the braided
tensor product algebra $B\und\tens B$ (the underscore is to distinguish it from
the usual tensor product algebra) and that $S_B$ is braided-antimultiplicative
(deduced from Proposition~2),
\ceqn{copmult}{\Delta_B(bc)=\sum b\Bo c\Bo\bo\tens b\Bt\bo c\Bt
\CR(b\Bt\bt\tens c\Bo\bt)\\
S_B(bc)=\sum (S_B c\bo)(S_B b\bo)\CR(b\bt\tens c\bt)}
in the case of $B\in \CM^H$, and similarly
\ceqn{lcopmult}{\Delta_B(bc)=\sum \CR(c\Bo\bo\tens b\Bt\bo)b\Bo c\Bo\bt\tens
b\Bt\bt c\Bt\\
S_B(bc)=\sum \CR(c\bo\tens b\bo)(S_B c\bt)(S_B b\bt)}
in the case of ${}^H\CM$. Here we would like to clarify
\cite[eqn.~(1.9)]{FisMon:sch}: the property for $\Delta_B$ is part of the
definition while the property for $S_B$ is derived from Proposition~2 due to
the author. As well as studying such braided groups, the paper \cite{Ma:bg}
introduced a `transmutation' procedure $B(\ ,\ )$ for their construction by
means of a generalised Tannaka-Krein reconstruction theorem. The property
(\ref{copmult}) is verified directly in \cite[Prop.~A.5]{Ma:bg} for the example
$B(H,H)\in \CM^H$ associated canonically to $H$ itself.

\new{4. \sl Dualisation.} We recall that it is one of the standard ideas behind
the theory of ordinary  Hopf algebras that for every construction in the
category of $k$-modules based on hypotheses and proofs by commutative diagrams,
there is a dual one by reversing all arrows. The axioms of a Hopf algebra are
self-dual in this respect, while the axioms of a comodule are dual to those of
a module, the construction of smash products dual to that of smash coproducts,
etc. For constructions which can be expressed as commuting diagrams it is
therefore redundant to give both module and comodule versions and one could not
really be considered mathematically new once the other is known. The
constructions to be considered here are all clearly of this diagrammatic type.
Note that it is {\em not any specific} Hopf algebra which is being dualised
here. This is analogous to the convention that a result for left modules is not
new when reworked for right modules, etc., the symmetry in this case being
reversal of $\tens$ in the category of $k$-modules.

\new{5. \sl Bosonisation.} One of the main results to date about braided groups
is the  `bosonisation theorem' introduced by the author in \cite{Ma:bos}. This
generalises the Jordan-Wigner bosonisation transform for $\Z_2$-graded systems
in physics, and associates to every Hopf algebra $B$ in the braided category of
representations of $H$ an equivalent ordinary Hopf algebra $B\lbiprod H$. The
paper \cite{Ma:bos} focused on the case where $B$ is in ${}_H\CM$ where $H$ is
quasitriangular.

\begin{propos}\cite[Thm~4.1]{Ma:bos} Let $H$ be quasitriangular. If $B\in
{}_H\CM$ is a braided group then $B\lbiprod H$ defined as the smash product by
the canonical action $\la$ of $H$ (by which $B$ is an object) and coproduct
$\Delta b=\sum b\Bo\CR\ut\tens \CR\uo\la b\Bt$ for all $b\in B$ is a Hopf
algebra over $k$. Here  $\CR=\sum\CR\uo\tens\CR\ut\in H\tens H$ is the
quasitriangular structure.
\end{propos}

As emphasised later in \cite[Thm~4.2]{Ma:skl}, the coproduct is also a smash
coproduct, namely by the {\em induced coaction} $\beta(b)=\CR\ut\tens\CR\uo\la
b$ introduced by the author in \cite[Thm~3.1]{Ma:dou}. Hopf algebras which are
both smash products and smash coproducts by the same (co)acting Hopf algebra
have been called `biproducts' by Radford\cite{Rad:str}. The observation that
bosonisations $B\lbiprod H$ are examples of this type is due to the author in
\cite{Ma:skl}. This is not quite clear from \cite[Remark~1.17]{FisMon:sch}.

If $H$ is a dual quasitriangular Hopf algebra then it is a trivial matter to
dualise the above construction in the sense of Section~4 above. If $B$ is a
braided group in $\CM^H$ we obviously make a smash coproduct $H\rbiprod B$ by
the canonical coaction of $H$ and smash product by the {\em induced (right)
action} $b\ra h=\sum b\bo\CR(b\bt\tens h)$. Likewise if $B$ in ${}^H\CM$ as in
\cite{FisMon:sch} we make a left handed smash coproduct $B\lbiprod H$ by the
canonical coaction and smash  product by the induced left action $h\la b=\sum
\CR(b\bo\tens h)b\bt$.
These are clearly the variants of the bosonisation theorem for Hopf algebras in
$\CM^H$ and ${}^H\CM$. Their further properties along the lines of
\cite{Ma:bos} are obtained by turning the diagrammatic proofs of these
properties up-side-down or reversing arrows. For example:

\begin{propos}cf\cite[Thm~4.2]{Ma:bos} Let $H$ be a dual quasitriangular Hopf
algebra and $B$ a Hopf algebra in $\CM^H$. The $B$-comodules in $\CM^H$ can be
identified canonically with $H\rbiprod B$-comodules as monoidal categories over
${}_k\CM$.
\end{propos}
\proof Cf\cite{Ma:bos} a $B$-comodule in the category means a vector space
which is an $H$-comodule and a $B$-comodule which intertwines the $H$-coaction.
The corresponding coaction of $H\rbiprod B$ consists of the $B$-coaction
followed by the $H$-coaction. Using the properties of a dual quasitriangular
structure one finds that this is an identification of monoidal categories (i.e.
that the tensor product of comodules is respected). \endproof

Likewise, left comodules of $B\in{}^H\CM$ in the category can be identified
with left comodules of $B\lbiprod H$. The more conceptual proof of
Proposition~4 along the lines in \cite{Ma:bos} is that
\eqn{cobos}{B(H\rbiprod B,H)\isom B(H,H)\rcocross B}
where $B(\ ,\ )$ is the transmutation procedure  from \cite[Thm~4.2]{Ma:bg}
mentioned in Section~3, and the right hand side is the categorical smash
coproduct of Hopf algebras in $\CM^H$ with braided tensor product algebra
structure.  The corresponding categorical smash products were introduced in
\cite[Thm~2.4]{Ma:bos} for just this purpose in the module setting.

We consider these results to be variants of a single bosonisation theorem and
did not explicitly elaborate all (four) natural cases in \cite{Ma:bos}. This
differs from the presentation in \cite[Remark~1.17]{FisMon:sch} where the
historical order also appears to us unclear (from the published left module
version of the bosonisation theorem one at once derives the right module and
the comodule versions).
In any case, the first draft of \cite{FisMon:sch} shown to the author at the
LMS meeting in July 1992 did not contain exactly the left or right comodule
bosonisation formulae above,  but a version appropriate to a more unnatural
set-up involving $H^{\rm cop}$-modules, see \cite[Remark~1.16]{FisMon:sch}.  It
was explained to Susan Montgomery at the time that the correct dualisation of
\cite{Ma:bos} did not need such a formulation with $H^{\rm cop}$. We would like
to stress that this is as true for the right comodule bosonisation and
biproducts as for the left comodule bosonisation and biproducts. Otherwise
\cite[Remark~1.16]{FisMon:sch} could lead to some confusion.

\new{6. \sl Duality and bosonisation.} When $B,H$ are finite-dimensional then
not only are the module and comodule bosonisation constructions dual, the
specific examples are dual as ordinary Hopf algebras as well. One has
\cite[Lem~4.3]{Ma:mec}
\eqn{dualbos}{ H^*\rbiprod B^\star\isom(B\lbiprod H)^*}
where $B^\star=B^{*\rm op/cop}$ is the natural dual of $B\in {}_H\CM$ where $H$
is quasitriangular. This is the new algebraic result in \cite{Ma:mec} (which
otherwise develops an application).
Note that the most natural identification of categories is ${}_H\CM=\CM^{H^*}$
by evaluation of a right coaction to give a left action of the dual. This is
why we use the right comodule bosonisation in (\ref{dualbos}) and
\cite{Ma:mec}.

\new{7. \sl Biproducts.}  In case it is not clear from \cite{FisMon:sch} we
would like to stress that the picture for general biproducts $B\lbiprod H$
(denoted $B\star H$ in \cite[Prop.~1.15]{FisMon:sch}) in terms of Hopf algebras
in braided categories is due to the author in
\cite[Prop.~A.2]{Ma:skl}\cite{Ma:dou}. We no longer assume that $H$ is
quasitriangular or dual quasitriangular etc. Radford\cite{Rad:str} elaborated
the conditions for a smash product algebra and smash coproduct coalgebra
$B\lbiprod H$ to be a Hopf algebra, and characterised the Hopf algebras
obtained in this way as those equipped with a Hopf algebra projection.

Let $H$ be a Hopf algebra with bijective antipode. The braided category of left
crossed $H$-modules ${}^H_H\CM$ consists of objects $V$ which are both modules
and comodules compatible in a natural way
\eqn{lcrossmod}{ \sum h\o v\bo\tens h\t\la v\bt=\sum (h\o\la v)\bo h\t\tens
(h\o\la v)\bt.}
These were studied by Yetter in \cite{Yet:rep} as a generalisation of
Whitehead's crossed $G$-modules\cite{Whi:com} and shortly thereafter by the
author in \cite{Ma:dou}, where we observed that this category is nothing other
than a reformulation appropriate to the case of $H$ infinite-dimensional of the
well-known braided category ${}_{D(H)}\CM$, where $D(H)$ is the quantum double
construction of V.G. Drinfeld\cite{Dri}. One just views a module of $H^{*\rm
op}\subset D(H)$ as an $H$-comodule in the usual way.

Now consider $H$ acting and coacting on an algebra and coalgebra $B$ as in the
setting of \cite{Rad:str}. It was observed in \cite{Ma:dou} (below Cor.~2.3
there) that $B\in {}^H_H\CM$ according to (\ref{lcrossmod}) is precisely one of
Radford's principal conditions for $B\lbiprod H$ to be a bialgebra. The
supplementary conditions that $B$ is an $H$-module coalgebra and an
$H$-comodule algebra in \cite{Rad:str} were omitted from \cite[Cor~2.3]{Ma:dou}
but should be understood. In \cite{Ma:skl} (see proof of Proposition~A.2) we
noted further that the second of Radford's principal conditions also has a
categorical meaning, namely that $B$ is a bialgebra in ${}^H_H\CM$. Radford's
condition for an antipode corresponds to $B$ a Hopf algebra in ${}^H_H\CM$.
Hence

\begin{propos}\cite[Prop.~A.2]{Ma:skl} Let $H$ be a Hopf algebra with bijective
antipode. If $B$ is a Hopf algebra in the braided category ${}^H_H\CM$ then
$B\lbiprod H$ is an ordinary Hopf algebra with projection $B\lbiprod H\to H$.
Conversely, every Hopf algebra with projection to $H$ is of the form $B\lbiprod
H$ for $B$ a Hopf algebra in ${}^H_H\CM$.
\end{propos}

As in \cite{Ma:dou}, we emphasised the ${}_{D(H)}\CM$ point of view but also
explained that the same result holds for ${}^H_H\CM$ in the
infinite-dimensional case. The work \cite{Ma:skl} was circulated in
January~1992 and its original form is archived on ftp.kurims.kyoto-u.ac.jp as
kyoto-net/92-02-07-majid. The right handed version for $B\in \CM^H_H$ works
just as well and gives an ordinary Hopf algebra $H\rbiprod B$ by right handed
smash product and coproduct. The right-handed version of (\ref{lcrossmod}) is
\eqn{rcrossmod}{\sum v\bo\ra h\o\tens v\bt h\t=\sum (v\ra h\t)\bo\tens h\o(v\ra
h\t )\bt}
and when $H$ is finite-dimensional one has ${}_{H^*}^{H^*}\CM=\CM^H_H$ by the
standard identifications.

\new{8. \sl Relating bosonisation and biproducts.} The biproduct point of view
on our bosonisation constructions was  emphasised in \cite{FisMon:sch} based on
an identification ${}^H\CM={}^H_H\CM$ which is stated several times
\cite[p.594, eqn~(1.16), below Prop~1.15, Remark~1.16]{FisMon:sch}.
Indeed, the introduction of this paper does not discuss the author's
bosonisation work at all on this basis. Unfortunately, such an   identification
is incorrect and rather misleading because bosonisation and biproducts are not
the same that one can eliminate one and replace it by the other. This is the
principal confusion in \cite{FisMon:sch} which we clarify in this note. It is
important because bosonisations have many remarkable properties (such as
Proposition~4 above) which are not true for general biproducts.

The functor which connects bosonisation to biproducts is in any case due to the
author. In \cite[Prop.~3.1]{Ma:dou} we showed that if $H$ is quasitriangular
then there is a functor of braided monoidal categories
\eqn{func}{{}_H\CM\hookrightarrow {}^H_H\CM,\quad (V,\la)\mapsto
(V,\la,\CR(\la)),\quad \CR(\la)(v)=\sum\CR\ut\tens\CR\uo\la v.}
This adds to an action $\la$ the induced coaction as used in Proposition~3. One
verifies that the two fit together to form a crossed module, and that this
identification respects tensor products. This means that an algebra, coalgebra,
Hopf algebra etc in ${}_H\CM$ can be viewed in ${}^H_H\CM$. It is clear that
$B\lbiprod H$ constructed by bosonisation can also be viewed as an example of a
biproduct as already noted above and in \cite{Ma:skl}. This functor in
\cite{Ma:dou} was the first of its kind. As an application we showed that the
quantum double $D(H)$ was itself a biproduct (in fact, a bosonisation).

The version when $H$ is dual quasitriangular is obviously the functor
\eqn{cofunc}{\CM^H\hookrightarrow \CM^H_H, \quad (V,\beta)\mapsto
(V,\beta,\la),\quad v\ra h=v\bo \CR(v\bt\tens h)}
for the right comodule version, or
\eqn{lcofunc}{{}^H\CM\hookrightarrow {}^H_H\CM, \quad (V,\beta)\mapsto
(V,\beta,\la),\quad h\la v=\CR(v\bo\tens h)v\bt}
for the left comodule version. Here $\la$ is the induced action as in Section~5
and one can check directly from (\ref{dqua}) that this identification respects
tensor products. This is the variant or our result from \cite{Ma:dou} which was
used in the final version of  \cite{FisMon:sch}. We note in passing that the
category ${}^H_H\CM$ is an example of a more general construction of the
`Pontryagin dual'\cite{Ma:rep} or `double' (in the sense of V.G. Drinfeld) of
any monoidal category. The above functors likewise have a categorical origin
\cite[Prop~2.5]{Ma:cat}.

\begin{propos} Let $H$ be a quasitriangular Hopf algebra with $H\ne k$. Then
the functor ${}_H\CM\to {}^H_H\CM$ introduced in \cite{Ma:dou} is never an
isomorphism. Likewise, let $H$ be a dual quasitriangular Hopf algebra with
$H\ne k$. Then ${}^H\CM\to {}^H_H\CM$ is never an isomorphism.
\end{propos}
\proof This is clear in the finite-dimensional case from the construction in
\cite{Ma:dou}, where this functor was introduced as pull back along a Hopf
algebra projection $D(H)\to H$. Since $D(H)$ as a vector space is $H^*\tens H$,
this can never be isomorphism. An isomorphism of categories would, by
Tannaka-Krein reconstruction, require such an isomorphism. This is the
conceptual reason. For a formal proof which includes the infinite-dimensional
case, consider $H\in {}^H_H\CM$ by the left regular coaction $\Delta$ and left
adjoint action. If in the image of the first functor (with $H$
quasitriangular), then $\sum h\o\tens h\t=\sum \CR\ut \tens \CR\uo\o h
S\CR\uo\t$ for all $h$ in $H$. Applying $\eps$ to the second factor tells us
that $h=\eps(h)$ for all $h$, i.e. $H=k$. This object in ${}^H_H\CM$ {\em can}
be in the image of the second functor, but this is iff the dual quasitriangular
structure is trivial and $H$ commutative. On other hand, consider $H\in
{}^H_H\CM$ by the left regular action and left adjoint coaction. If in the
image of the second functor (with $H$ dual quasitriangular) then $hg= \CR(g\o
Sg\th\tens h)g\t$. Setting $g=1$ tells us that $h=\eps(h)$ again, hence $H=k$.
This object can be in the image of the first functor, but this is iff the
quasitriangular structure is trivial and $H$ cocommutative. \endproof

This means in turn that general `biproducts' associated to $B\in {}^H_H\CM$ are
much more general than the Hopf algebras obtained by bosonisation when $B\in
{}_H\CM$ for $H$ quasitriangular or $B\in {}^H\CM$ for $H$ dual
quasitriangular, both of the latter being viewable as examples of biproducts.

\new{9. \sl Symmetric braiding.} We would like to stress that inspite of the
term `symmetric braiding' used in \cite{FisMon:sch} for the CT structure of a
cotriangular (CT) Hopf algebra, there is no canonical braid group action in
this setting. This can be misleading. For example, the colour-Lie algebras etc.
and their universal enveloping algebra studied in \cite{FisMon:sch} are defined
in the obvious way with transposition replaced by a representation of the
symmetric group. There {\em is} a theory of truly braided-Lie algebras and
their enveloping algebras introduced in \cite{Ma:lie} but the definition of the
Jacobi identity is rather more complicated in the braided case.

\new{10. \sl Reconstruction.} In the 1990 paper \cite[Thm~2.2]{Ma:bg} we
introduced and proved a very general reconstruction theorem which yields a Hopf
algebra $\aut(\CC,\omega,\CV)\in\CV$ associated to a monoidal functor
$\omega:\CC\to\CV$ between a monoidal category $\CC$ and a braided monoidal
category $\CV$. This is the {\em automorphism (or endomorphism) braided group}
of a functor and is due to the author as a significant generalisation of usual
Tannaka-Krein ideas. It is constructed in \cite{Ma:bg} as representing object
for $\Nat(\omega,\omega(\ ))$. As an application of it we obtained in
\cite{Ma:bg} the transmutation construction $B(\ ,H)$ mentioned above, which
turns any ordinary  Hopf algebra mapping to $H$ into a braided group in
$\CM^H$. Here $H$ is dual quasitriangular. So far, this transmutation remains
one of the main constructions for Hopf algebras in braided categories such as
$\CM^H$.

Recently there appeared an interesting preprint \cite{Par:rec} in the
introduction of which the following question is posed: Let $H$ be dual
quasitriangular, $B$ a Hopf algebra in $\CM^H$ and $\omega$ the forgetful
functor from $B$-comodules in $\CM^H$ to $\CM^H$. What is the automorphism
braided group $\aut(\omega)$ in $\CM^H$ reconstructed in this case? Our
generalised reconstruction  work is recalled explicitly in
\cite[Sec.~3.4.1]{Par:rec} and a partial answer (for the coalgebra structure
when $B$ is only a coalgebra) is obtained\cite[Cor~5.7]{Par:rec} as a main
result of the paper.

Here we point out that this question was already answered (the full braided
group structure of $\aut(\omega)$) in the course of our 1991 paper on
bosonisation\cite{Ma:bos}. Some refinements of the problem to the case of
`limited reconstruction' over a control category in \cite{Par:rec} remain
interesting and will not be addressed here.

\begin{propos} Let $H$ be dual quasitriangular and $B$ a Hopf algebra in
$\CM^H$. Then the forgetful functor $\omega$ from $B$-comodules in $\CM^H$ to
$\CM^H$ has as automorphism braided group the Hopf algebra $B(H,H)\rcocross B$
in $\CM^H$ mentioned in (\ref{cobos}). It has the smash coproduct coalgebra and
braided tensor product algebra, and is a transmutation of the bosonisation
$H\rbiprod B$ of $B$.
\end{propos}
\proof  Under the equivalence in Proposition~4, the forgetful functor $\omega$
becomes the functor induced by push-out along the canonical Hopf algebra map
$H\rbiprod B\to H$ defined by the counit of $B$. But the automorphism braided
group of a functor induced by push out is exactly the definition of the
transmutation construction $B(\ ,H)$. So the answer is exactly the
transmutation $B(H\rbiprod B,H)$. But the abstract definition of bosonisation
$H\rbiprod B$ is (\ref{cobos}) i.e., exactly such that its transmutation is
$B(H,H)\rcocross B$. These are exactly the conceptual steps (in comodule form)
which led to the author's bosonisation theory \cite{Ma:bos} in the first place.
\endproof

This demonstrates how one may use the bosonisation theory of \cite{Ma:bos}: we
convert our problem for the braided group $B$ to one for its equivalent
ordinary Hopf algebra $H\rbiprod B$. Explicitly, the braided group $B(H,H)$
associated to $H$ was introduced as the automorphism braided group of the
identity functor from $\CM^H$ to itself\cite{Ma:eul}\cite{Ma:bg} and
corresponds to $B=k$. Its structure is $H$ as a coalgebra, with the right
adjoint coaction and  modified product\cite{Ma:bg}
\eqn{BHH}{ \sum h\bo\tens h\bt=\sum h\t \tens (Sh\o)h\th,\quad h\cdot g=\sum
h\t g\t \CR((Sh\o)h\th\tens Sg\o)}
in terms of the structure of $H$. This $B(H,H)$ coacts on any $B$ by the same
map $\beta$ by which $H$ coacts on $B$ as an object (the {\em tautological
coaction}). The fact that one can then make a (braided) smash coproduct by this
 and still obtain a Hopf algebra in the braided category with the braided
tensor product algebra structure reflects the fact that $B(H,H)$ is
braided-commutative with respect to $B$ in a certain (unobvious) sense
introduced in\cite{Ma:bg}. This was the key idea behind the construction in
\cite[Sec.~2]{Ma:bos}. This $B(H,H)\rcocross B$ has product and coproduct
defined diagrammatically cf.\cite[Sec.~2]{Ma:bos}
\eqn{rcross}{\epsfbox{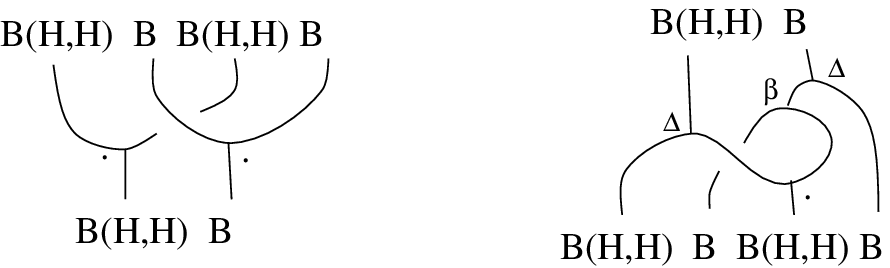}}
where $\Psi=\epsfbox{braid.eps}$ is the braiding. The notation is explained
further in \cite{Ma:bos}.
In our particular case in the category $\CM^H$ it means
\alignn{rec}{ (h\tens b)(g\tens c)\equad &&=\sum h\cdot g\bo \tens b\bo
c\CR(b\bt\tens g\bt)\nonumber \\
&&=\sum h\t g\th \tens b\bo c \CR((Sh\o)h\th\tens  Sg\t)\CR(b\bt\tens
(Sg\o)g\fo) \\
\Delta (h\tens b)\equad &&=\sum h\o\tens  b\Bo\bo\bo\tens h\t\bo \cdot
b\Bo\bt\tens b\Bt\CR(h\t\bt\tens b\Bo\bo\bt)\nonumber\\
&&=\sum h\o\tens b\Bo\bo\tens h\t b\Bo\bt\tens b\Bt}
where we evaluated further in terms of $B,H$. The counit is the tensor product
one and there is an antipode as well. The coproduct comes out just the same as
for the bosonisation $H\rbiprod B$ (the usual smash coproduct by the coaction
of $H$ on $B$ as an object in $\CM^H$) because the transmutation procedure $B(\
,H)$ does not change the coalgebra. This is just the dual of the calculation of
$B\lcross B(H,H)$ in \cite[Thm~3.2]{Ma:bos} for the module version.

We recall that the diagrammatic smash products and coproducts where algebraic
information `flows' along braids were introduced by the author in
\cite{Ma:bos}. The adjoint action for braided groups was introduced and studied
by the author in \cite[Prop.~3.1]{Ma:lie} as a foundation for a theory of
braided Lie algebras. We proved that it makes $B$ a braided module algebra in
the category, etc. We also introduced the adjoint coaction (the adjoint action
turned up-side-down) in \cite[Appendix]{Ma:lin}. Some of these constructions
are used \cite[Sec~2.5]{Par:rec} in a slightly generalised form following the
same proofs. We would like to stress also that this novel diagrammatic way of
working is not a trivial generalisation of  `wiring diagrams' for usual linear
algebra in $Vec$ (as in Penrose's spin networks and \cite{Yet:rep}) because
under and over crossings now represent a nontrivial operator $\Psi$; we add the
rules for this based on functoriality and the coherence theorem for symmetric
or braided categories\cite{JoyStr:bra}.
See \cite{Ma:bos}.

\new{11. \sl Example.} A trivial example of Proposition~7 is to the case of
reconstruction of a super-Hopf algebra $\Z_2'\rcocross B$ from the category of
super $B$-modules and its forgetful functor. Here $\Z_2'$ is the dual of the
triangular Hopf algebra introduced by the author in
\cite[Prop~6.1]{Ma:exa}\cite[Ex.~1.1]{Ma:tra} as generator of the category
${\rm SuperVec}$ of superspaces with their $\Z_2$-graded transposition. This
novel application of quasitriangular Hopf algebras was generalised in 1991 in
\cite{Ma:any} to generate the braided category of anyonic or $\Z_n$-graded
vector spaces introduced there. Unfortunately, in all these examples the
adjoint coaction of $H$ is trivial and $B(H,H)=H$ is viewed trivially in
$\CM^H$. Hence the algebra structure of $B(H,H)\rcocross B$ is the usual tensor
product one (and its smash coproduct the usual one as well). The result is a
Hopf algebra in $\CM^H$ just because $B$ is. The same applies for Hopf algebras
in the braided category of $\Z$-graded vector spaces.

To give a more non-trivial example, let $q\in k^*$ and $H=GL_q(2)$ defined as
$k\<\alpha,\beta,\gamma,\delta,C^{-1}\>$ modulo the relations
\ceqn{qmat}{\alpha \beta=q^{-1}\beta \alpha,\quad \alpha \gamma=q^{-1}\gamma
\alpha,\quad \beta\delta=q^{-1}\delta \beta,\quad
\gamma\delta=q^{-1}\delta\gamma\\
\beta\gamma=\gamma\beta,\quad
\alpha\delta-\delta\alpha=(q^{-1}-q)\beta\gamma,\quad
C=\alpha\delta-q^{-1}\beta\gamma}
essentially as in \cite{Dri}\cite{FRT:lie} for $SU_q(2)$. We equip it now with
dual quasitriangular structure determined by the associated solution $R$ of the
quantum Yang-Baxter equations. More precisely (for our application) we take $R$
with a non-standard normalisation as explained in\cite{Ma:lin}, fixed instead
by $\CR(C\tens C)=q^6$. Note that one cannot set $C=1$ as one would for the
usual dual quasitriangular Hopf algebra $SU_q(2)$.

The braided group $B(GL_q(2),GL_q(2))=BGL_q(2)$ is likewise a variant of the
braided group $BSU_q(2)$ introduced by the author in
\cite{Ma:eul}\cite{Ma:exa}. We define it as $k\<a,b,c,d,D^{-1}\>$ modulo the
relations
\ceqn{bmat}{ba=q^2ab,\quad ca=q^{-2}ac,\quad d a=ad,\qquad
bc=cb+(1-q^{-2})a(d-a)\\
d b=bd+(1-q^{-2})ab,\quad cd=d c+(1-q^{-2})ca,\quad D=ad-q^2cb}
It has `matrix' coproduct $\Delta\vecu=\vecu\tens\vecu$ and $\Delta D=D\tens D$
when we regard the generators as a matrix $\vecu=\pmatrix{a&b\cr c&d}$. The
braided group antipode for $\vecu$ is as for $BSU_q(2)$ in
\cite{Ma:eul}\cite{Ma:exa} times $D^{-1}$. The braiding $\Psi$ between the
generators is also as listed for $BSU_q(2)$ in \cite{Ma:eul}\cite{Ma:exa}. This
$BGL_q(2)$ lives in $\CM^{GL_q(2)}$ with a coaction which has the same `matrix
conjugation' form on the generators $\vecu$ as the right adjoint coaction of
$GL_q(2)$.

Let $B={\Bbb A}_q^2=k\<x,y\>/(yx-qxy)$ the $q$-deformed plane with right
coaction of $GL_q(2)$ given by transformation of the $(x,y)$ as a row vector by
the $GL_q(2)$ generators as a matrix, i.e. $\beta(x)=x\tens \alpha+y\tens
\gamma$ and $\beta(y)=x\tens\beta+y\tens\delta$. One of the first applications
of braided groups to physics was to show that this `quantum-braided plane'
${\Bbb A}_q^2$ is  a Hopf algebra in $\CM^{GL_q(2)}$ with linear `coaddition'
\cite{Ma:poi}
\ceqn{qplane}{ \Psi(x\tens x)=q^2 x\tens x,\quad \Psi(x\tens y)=q y\tens x\\
\Psi(y\tens y)=q^2 y\tens y,\quad \Psi(y\tens x)=q x\tens y+(q^2-1)y\tens x\\
\Delta x=x\tens 1+1\tens x,\quad \Delta y=y\tens 1+1\tens y,\quad \eps x=0=\eps
y,\quad Sx=-x,\quad Sy=-y.}
This result is due to the author in \cite{Ma:poi}, where $GL_q(2)$ above is
formulated as $\widetilde{SU_q(2)}$, the `dilatonic' central extension.

We use the same matrix transformation for the braided coaction of $BGL_q(2)$ on
${\Bbb A}_q^2$. Under this, ${\Bbb A}_q^2$ becomes a right comodule algebra in
the braided category\cite[Prop.~3.7]{Ma:lin}.

\begin{example} The automorphism braided group $BGL_q(2)\rcocross {\Bbb A}_q^2$
in $\CM^{GL_q(2)}$ is generated by $BGL_q(2)$ and the quantum-braided plane
${\Bbb A}_q^2$ as subalgebras with the cross relations
\cmath{xa=ax,\quad ya=bx(q-q^{-1}) + ay,\quad xb=q^{-1}bx,\quad yb=qby\\
xc=qcx,\quad yc=(1-q^{-2})(d-a)x + q^{-1}cy ,\quad xd=dx,\quad
yd=dy-q^{-2}(q-q^{-1})bx}
It has the matrix coproduct of $BGL_q(2)$ and
\[ \Delta x=x\tens a+y\tens c+1\tens x,\quad \Delta y=x\tens b+y\tens d+1\tens
y\]
extended as a braided group in $\CM^{GL_q(2)}$.
\end{example}
\proof The cross relations are exactly the braided tensor product algebra as in
Section~2, computed for the present setting  in terms of $R$ in
\cite[Lem.~3.4]{Ma:lin}. This gives the relations shown. For the coproduct we
know that we have the same form as the smash coproduct by the coaction of
$GL_q(2)$ on ${\Bbb A}_q^2$ but viewed now as a coaction of $BGL_q(2)$. To
extend the coproduct to products of the generators we use its
braided-multiplicativity, with $\Psi$ determined from the coaction. This was
computed in terms of $R$ in \cite[Prop.~3.2]{Ma:lin} and in our case is
\ceqn{psiax}{\Psi(a\tens x)=x\tens a+(1-q^2)y\tens c,\  \Psi(b\tens
x)=q^{-1}x\tens b+(q-q^{-1})y\tens (a-d),\  \Psi(c\tens x)=qx\tens c\\
\Psi(d\tens x)=x\tens d+ (1-q^{-2})y\tens c,\quad
\Psi(\pmatrix{a &b\cr c&d}\tens y)=y\tens \pmatrix{a&qb\cr q^{-1}c&d}}
while the braiding $\Psi(x\tens a)=a\tens x$ etc., has just the same form as
the cross relations already given. It is enough to specify the coproduct and
braiding on the generators since the braiding $\Psi$ itself extends
`multiplicatively' by functoriality and the Hexagon coherence identities, as
explained in \cite{Ma:exa}. \endproof

The construction of linear braided groups such as the quantum-braided plane
works for general quantum planes associated to suitable matrix
data\cite{Ma:poi}. Another example is the 1-dimensional case $B={\Bbb
A}_q=k[x]$, the braided line\cite{Ma:sol}. Such `linear braided groups'  have
been very extensively studied since \cite{Ma:poi} as the true foundation for
$q$-deformed geometry. See \cite{Ma:varen} for a review. Their bosonisations
were used in \cite{Ma:poi} to define inhomogeneous quantum groups and are also
extensively studied since then.

\baselineskip 20pt


\begin{thebibliography}{Maj94d}

\bibitem[D1]{Dri}
V.G. Drinfeld.
\newblock Quantum groups.
\newblock In A.~Gleason, editor, {\em Proceedings of the {ICM}}, pages
  798--820, Rhode Island, 1987. AMS.

\bibitem[D2]{Dri:alm}
V.G. Drinfeld.
\newblock On almost-cocommutative {H}opf algebras.
\newblock {\em Leningrad Math. J.}, 1:321--342, 1990.

\bibitem[FM]{FisMon:sch}
D.~Fischman and S.~Montgomery.
\newblock A {S}chur double centralizer theorem for cotriangular {H}opf algebras
  and generalized {L}ie algebras.
\newblock {\em J. Algebra}, 168:594--614, 1994. Recvd June 1992 (revsd later).


\bibitem[FRT]{FRT:lie}
L.D. Faddeev, N.Yu. Reshetikhin, and L.A. Takhtajan.
\newblock Quantization of {L}ie groups and {L}ie algebras.
\newblock {\em Leningrad Math. J.}, 1:193--225, 1990.

\bibitem[JS]{JoyStr:bra}
A.~Joyal and R.~Street.
\newblock Braided monoidal categories.
\newblock Mathematics Reports 86008, Macquarie University, 1986.

\bibitem[LT]{LarTow:two}
R.G. Larson and J.~Towber.
\newblock Two dual classes of bialgebras related to the concepts of `quantum
  group' and `quantum {L}ie algebra'.
\newblock {\em Commun. Algebra}, 19(12):3295--3345, 1991.

\bibitem[LM]{LyuMa:bra}
V.~Lyubashenko and~S. Majid.
\newblock Braided groups and quantum Fourier transform.
\newblock {\em J. Algebra}, 166:506--528, 1994. Recvd December 1991.

\bibitem[M]{Ma:introm}
S.~Majid.
\newblock Algebras and {H}opf algebras in braided categories.
\newblock volume 158 of {\em Lec. Notes in Pure and Appl. Math}, pages 55--105.
  Marcel Dekker, 1994.


\bibitem[M1]{Ma:seq}
S.~Majid.
\newblock More examples of bicrossproduct and double cross product {H}opf
  algebras.
\newblock {\em Isr. J. Math}, 72:133--148, 1990. Recvd July 1989.

\bibitem[M2]{Ma:qua}
S.~Majid.
\newblock Quasitriangular Hopf algebras and Yang-Baxter equations.
\newblock {\em Int. J. Mod. Phys.}, 5:1--91, 1990. Recvd April 1989.

\bibitem[M3]{Ma:pro}
S.~Majid.
\newblock Quantum groups and quantum probability.
\newblock In {\em Quantum Probability and Related Topics VI (Proc. Trento,
  1989)}, pages 333--358. World Scientific.

\bibitem[M4]{Ma:tan}
S.~Majid.
\newblock Tannaka-{K}rein theorem for quasi{H}opf algebras and other results.
\newblock {\em Contemp. Maths}, 134:219--232, 1992. Recvd September 1990.

\bibitem[M5]{Ma:eul}
S.~Majid.
\newblock Rank of quantum groups and braided groups in dual form.
\newblock In {\em Proceedings Quantum Groups, St Petersburg, September 1990};
  {\em Lec. Notes. in Math.}, 1510:79--89, 1992. Springer.

\bibitem[M6]{Ma:bg}
S.~Majid.
\newblock Braided groups.
\newblock {\em J. Pure and Applied Algebra}, 86:187--221, 1993. Recvd March
1991;
Presented at the AMS Biannual Meeting, San Francisco, January 1991.

\bibitem[M7]{Ma:dou}
S.~Majid.
\newblock Doubles of quasitriangular {H}opf algebras.
\newblock {\em Commun. Algebra}, 19(11):3061--3073, 1991. Recvd February 1990.

\bibitem[M8]{Ma:skl}
S.~Majid.
\newblock Braided matrix structure of the {S}klyanin algebra and of the quantum
  {L}orentz group.
\newblock {\em Commun. Math. Phys.}, 156:607--638, 1993. Recvd February 1992.

\bibitem[M9]{Ma:bos}
S.~Majid.
\newblock Cross products by braided groups and bosonization.
\newblock {\em J. Algebra}, 163:165--190, 1994. Recvd June 1991.


\bibitem[M10]{Ma:tra}
S.~Majid.
\newblock Transmutation theory and rank for quantum braided groups.
\newblock {\em Math. Proc. Camb. Phil. Soc.}, 113:45--70, 1993. Recvd December
1991.

\bibitem[M11]{Ma:mec}
S.~Majid.
\newblock The quantum double as quantum mechanics.
\newblock {\em J. Geom. Phys.}, 13:169--202, 1994. Recvd October 1992.

\bibitem[M12]{Ma:rep}
S.~Majid.
\newblock Representations, duals and quantum doubles of monoidal categories.
\newblock In {\em Proceedings of the Srni Winter School, January 1990};
{\em Suppl. Rend. Circ. Mat. Palermo, Ser. II}, 26:197--206, 1991.

\bibitem[M13]{Ma:cat}
S.~Majid.
\newblock Braided groups and duals of monoidal categories.
\newblock In {\em Proc. Category Theory, Montreal, 1991};
{\em Canad. Math. Soc. Conf. Proc.}, 13:329--343, 1992.

\bibitem[M14]{Ma:lie}
S.~Majid.
\newblock Quantum and braided {L}ie algebras.
\newblock {\em J. Geom. Phys.}, 13:307--356, 1994. Recvd April 1993.

\bibitem[M15]{Ma:any}
S.~Majid.
\newblock Anyonic quantum groups.
\newblock In Z. Oziewicz et al., eds, {\em Spinors, Twistors, Clifford
Algebras and Quantum Deformations}, pages 327--336, 1993. Kluwer.

\bibitem[M16]{Ma:sol}
S.~Majid.
\newblock Solutions of the {Y}ang-{B}axter equations from braided-{L}ie
  algebras and braided groups, 1993.
\newblock To appear in {\em J. Knot Th. Ram.}


\bibitem[M17]{Ma:exa}
S.~Majid.
\newblock Examples of braided groups and braided matrices.
\newblock {\em J. Math. Phys.}, 32:3246--3253, 1991.

\bibitem[M18]{Ma:poi}
S.~Majid.
\newblock Braided momentum in the {$q$}-{P}oincar{\'e} group.
\newblock {\em J. Math. Phys.}, 34:2045--2058, 1993.

\bibitem[M19]{Ma:lin}
S.~Majid.
\newblock Quantum and braided linear algebra.
\newblock {\em J. Math. Phys.}, 34:1176--1196, 1993.

\bibitem[M20]{Ma:varen}
S.~Majid.
\newblock Introduction to braided geometry and {q-M}inkowski space.
\newblock To appear in {\em Proceedings of the School on Quantum Groups,
{V}arenna,
  {I}taly, {J}une 1994}. 60 pages.

\bibitem[P1]{Par:non}
B.~Pareigis.
\newblock A non-commutative non-cocommutative Hopf algebra in `Nature'.
\newblock {\em J. Algebra}, 70:356--374, 1981.


\bibitem[P2]{Par:rec}
B.~Pareigis.
\newblock Reconstruction of hidden symmetries.
\newblock {\em Preprint}, December, 1994.

\bibitem[R]{Rad:str}
D.~Radford.
\newblock The structure of {H}opf algebras with a projection.
\newblock {\em J. Algebra}, 92:322--347, 1985.

\bibitem[S]{Swe:hop}
M.E. Sweedler.
\newblock {\em {H}opf Algebras}.
\newblock Benjamin, 1969.

\bibitem[W]{Whi:com}
J.H.C. Whitehead.
\newblock Combinatorial homotopy, {II}.
\newblock {\em Bull. Amer. Math. Soc.}, 55:453--496, 1949.

\bibitem[Y]{Yet:rep}
D.N. Yetter.
\newblock Quantum groups and representations of monoidal categories.
\newblock {\em Math. Proc. Camb. Phil. Soc.}, 108:261--290, 1990.

\end{thebibliography}

\end{document}